\documentclass[a4paper]{article}%
\usepackage{amsfonts}
\usepackage{amsmath}
\usepackage{amssymb}
\usepackage{graphicx}%
\providecommand{\U}[1]{\protect\rule{.1in}{.1in}}
\newtheorem{theorem}{Theorem}

\newtheorem{corollary}[theorem]{Corollary}

\newtheorem{lemma}[theorem]{Lemma}

\newenvironment{proof}[1][Proof]{\noindent\textbf{#1.} }{\ \rule{0.5em}{0.5em}}
\begin{document}

\begin{center}
\textbf{A Fundamental Theorem of Calculus for }

\textbf{Second-order Directional Derivative}

Pisheng Ding

Illinois State University

\texttt{pding@ilstu.edu}
\end{center}

\begin{quote}
\textsc{Abstract} {\small Given a two-variable function \textit{f} without
critical points and a compact region \textit{R} bounded by two level curves of
\textit{f}, this note proves that the integral over \textit{R} of \textit{f}'s
second-order directional derivative in the tangential directions of the
interceding level curves is proportional to the rise in \textit{f}-value over
\textit{R}. Also discussed are variations on this result when critical points
are present or \textit{R} becomes unbounded. Several concrete examples
exemplify the theory.}
\end{quote}

\noindent Let $f$ be a real-valued $C^{2}$ function on an open connected
domain $\Omega\subset%
%TCIMACRO{\U{211d} }%
%BeginExpansion
\mathbb{R}
%EndExpansion
^{2}$. Suppose that $f$ has \textit{no} critical points and that $a<b$ are
values in $f[\Omega]$ such that $f^{-1}[[a,b]]$ is \textit{connected} and
\textit{compact} (in which case $f^{-1}[[a,b]]$ is diffeomorphic to
$f^{-1}[a]\times\lbrack a,b]$ and $f^{-1}[t]$ is a simple closed $C^{2}$ curve
for each $t\in\lbrack a,b]$). For $p\in\Omega$, let $\mathbf{T}(p)$ be
\textit{a} unit tangent to the level curve $f^{-1}[f(p)]$ at $p$. With
$D_{\mathbf{T}(p)}^{2}f(p)$ denoting the second-order directional derivative
of $f$ at $p$ in the direction $\mathbf{T}(p)$, our main result is the
following identity, unmistakably resembling the fundamental theorem of
calculus:%
\[
\iint\nolimits_{f^{-1}[[a,b]]}D_{\mathbf{T}}^{2}f\,dA=\pm2\pi(b-a)\,\text{,}%
\]
with the positive sign in effect iff $f^{-1}[b]$ encircles $f^{-1}[a]$. As we
shall see through several examples, the assumptions that $f^{-1}[[a,b]]$ be
connected and compact and that $f$ have no critical points can be relaxed,
allowing flexibility in application.

We establish this identity in \S 2 after treating some preparatory results in
\S 1. In \S 3, we show how this result can be adapted for a variety of situations.

\section{Second-order directional derivative and curvature of level curves}

We review a few key notions, aiming to conceptualize curvature of level curves
of a $C^{2}$ function $f$ in terms of its second-order directional derivative.

\subsection{Second-order directional derivative}

For notation, we often add a displacement vector $\mathbf{v}$ to an initial
point $p\in%
%TCIMACRO{\U{211d} }%
%BeginExpansion
\mathbb{R}
%EndExpansion
^{2}$ to express the terminal point $p+\mathbf{v}$.

For each $p\in\Omega$ and any unit vector $\mathbf{v}$, let
$D(s)=f(p+s\mathbf{v})$ for $s$ sufficiently small so that $p+s\mathbf{v}%
\in\Omega$. The first and second directional derivatives of $f$ at $p$ in the
direction $\mathbf{v}$, denoted by $D_{\mathbf{v}}f(p)$ and $D_{\mathbf{v}%
}^{2}f(p)$, are defined to be the two numbers $D^{\prime}(0)$ and
$D^{\prime\prime}(0)$. Using chain rule, we obtain the standard facts that%
\[
D_{\mathbf{v}}f(p)=\nabla f(p)\cdot\mathbf{v}\text{\quad and\quad
}D_{\mathbf{v}}^{2}f(p)=Q_{p}(\mathbf{v},\mathbf{v})\text{ ,}%
\]
where $Q_{p}$ is the quadratic form associated with the Hessian matrix%
\[
H_{p}=\left[
\begin{array}
[c]{cc}%
f_{xx}(p) & f_{xy}(p)\\
f_{yx}(p) & f_{yy}(p)
\end{array}
\right]  \text{.}%
\]
Simply put, $Q_{p}(\mathbf{v},\mathbf{v})=\mathbf{v}\cdot\left(
H_{p}\mathbf{v}\right)  $. We note another formula for $D_{\mathbf{v}}%
^{2}f(p)$.

\begin{lemma}
\label{Lemma 2nd Derivative/Gradient}Let $t\mapsto\mathbf{r}(t)$ be a curve in
$\Omega$ such that $\mathbf{r}(0)=p$ and $\mathbf{r}^{\prime}(0)=\mathbf{v}$.
Then,%
\[
D_{\mathbf{v}}^{2}f(p)=\left.  \frac{d}{dt}\right\vert _{t=0}\left(  \nabla
f(\mathbf{r}(t))\cdot\mathbf{v}\right)  =\left(  \left.  \frac{d}%
{dt}\right\vert _{t=0}\nabla f(\mathbf{r}(t)\right)  \cdot\mathbf{v\,}%
\text{.}
\]

\end{lemma}

\begin{proof}
The second equality is clear. Checking the first amounts to verifying that
$\left.  \frac{d}{dt}\right\vert _{t=0}\left(  \nabla f(\mathbf{r}%
(t))\cdot\mathbf{v}\right)  =Q_{p}(\mathbf{v},\mathbf{v})$ by chain rule.
\end{proof}

\subsection{Curvature of level curves}

Let $C$ denote a level curve of $f$, which is a regular $C^{2}$ curve (as, by
assumption, $f$ has no critical point). Install on $C$ the unit normal field
$\mathbf{N:}=-\nabla f/|\nabla f|$ and the unit tangent field $\mathbf{T:}%
=(-f_{y}\mathbf{e}_{1}+f_{x}\mathbf{e}_{2})/|\nabla f|$. (The frame
$(\mathbf{T},\mathbf{N})$ is positively-oriented.) The \textit{signed}%
\ \textit{curvature} $\kappa$ of $C$ at each point thereon is defined by the
equation $d\mathbf{T}/ds=\kappa\mathbf{N}$, where $s$ is arc length along $C $
with its increasing direction induced by $\mathbf{T}$. (The sign of $\kappa$
depends on the choice we make of $\mathbf{N}$, but not of $\mathbf{T}$.) For
$p\in\Omega$, let $\kappa(p)$ be the signed curvature of $f^{-1}[f(p)]$ at $p$.

\begin{lemma}
\label{Lemma Curvature of Level Curves}$\kappa(p)=D_{\mathbf{T}}%
^{2}f(p)/\left\vert \nabla f(p)\right\vert $.
\end{lemma}

\begin{proof}
Let $\gamma$ be the unit-speed parametrization of an arc on $C:=f^{-1}[f(p)]
$, with $\gamma(0)=p$ and $\gamma^{\prime}(0)=\mathbf{T}(p)$. By definition,
$\gamma^{\prime\prime}(0)=\kappa(p)\mathbf{N}(p)$. Because $C $ is a level
curve, the two vectors $\nabla f(\gamma(t))$ and $\gamma^{\prime}(t)$ are
always orthogonal; hence,%
\[
0=\frac{d}{dt}\left(  \nabla f(\gamma(t))\cdot\gamma^{\prime}(t)\right)
=\left(  \frac{d}{dt}\nabla f(\gamma(t))\right)  \cdot\gamma^{\prime
}(t)+\nabla f(\gamma(t))\cdot\gamma^{\prime\prime}(t)\text{\thinspace.}
\]
For the last two terms, note that, by Lemma
\ref{Lemma 2nd Derivative/Gradient},%
\[
\left(  \left.  \frac{d}{dt}\right\vert _{t=0}\nabla f(\gamma(t))\right)
\cdot\gamma^{\prime}(0)=D_{\mathbf{T}}^{2}f(p)\text{ ,}
\]
whereas%
\[
\nabla f(\gamma(0))\cdot\gamma^{\prime\prime}(0)=\nabla f(\gamma
(0))\cdot\left(  \kappa(p)\frac{-\nabla f(p)}{|\nabla f(p)|}\right)
=-\kappa(p)|\nabla f(p)|\,\text{.}
\]
The claimed formula now follows, since $|\nabla f(p)|\neq0$.
\end{proof}

In the literature, this formula for $\kappa$ is always written explicitly in
terms of the partial derivatives of $f$ and is typically derived (using the
implicit function theorem)\ from the curvature formula for graphs of
one-variable functions; see, e.g., \cite{Courant}. Not only do we formulate
the result in a conceptually simpler form, we have given a conceptually
simpler derivation not relying on \textit{any}\ formula for curvature other
than its definition.

\section{Integrating Second-order Directional Derivative}

We reiterate our assumptions, which will remain in effect in this section.

\smallskip

\noindent\textbf{Hypothesis. }$f$ is a $C^{2}$ function with \textit{no}
critical points on an open connected set $\Omega\subset%
%TCIMACRO{\U{211d} }%
%BeginExpansion
\mathbb{R}
%EndExpansion
^{2}$; $a<b$ are two numbers in $f[\Omega]$ such that $f^{-1}[[a,b]]$ is
\textit{connected} and \textit{compact}.

\smallskip

We first recall two facts and also introduce a notation.

\subsection{Two facts from calculus and geometry}

First, the gradient flow originating from $f^{-1}[a]$ induces a diffeomorphism
between $f^{-1}[[a,b]]$ and $f^{-1}[a]\times\lbrack a,b]$, allowing this
change of variables of integration:%
\begin{equation}
\iint\nolimits_{f^{-1}[[a,b]]}g\,dA=\int_{a}^{b}\left(  \int_{f^{-1}[t]}%
\frac{g}{|\nabla f|}ds\right)  dt \label{Eq Resolution}%
\end{equation}
where $ds$ is the arc length element along $f^{-1}[t]$; see
\cite[pp.\thinspace298--300]{Courant}.

Second, for a simple closed plane curve $C$ with signed curvature $\kappa$,
$\int_{C}\kappa ds=\pm2\pi$; see \cite[pp.\thinspace36-37]{do Carmo}. The sign
ambiguity is due to the dependence of $\kappa$ on orientation, with the
positive sign in force iff the chosen unit normal field on $C$ points inward
(relative to the Jordan domain enclosed by $C$). To encode the sign more
effectively, we introduce $\sigma:=\mathbf{n\cdot N}$, where $\mathbf{n}$ is
the inward unit normal field along $C$ and $\mathbf{N}$ is the chose normal
field. With $\sigma$ tracking orientation, $\int_{C}\kappa ds=2\pi\sigma$.

Turning to level curves of $f$, for $p\in f^{-1}[[a,b]]$, let $\mathbf{n}(p)$
be the inward unit normal at $p$ of the (simple closed) curve $f^{-1}[f(p)]$
and let $\sigma(p):=\mathbf{n}(p)\cdot\mathbf{N}(p)$. Being continuous and
integer-valued, $\sigma$ is constant on (the connected) $f^{-1}[[a,b]]$. (Note
that $\sigma\equiv1$ iff $f^{-1}[b]$ encloses $f^{-1}[a]$.) We then have, for
every $t\in\lbrack a,b]$,%
\begin{equation}
\int_{f^{-1}[t]}\kappa ds=2\pi\sigma\text{\thinspace.}
\label{Eq Integral-of-Curvature}%
\end{equation}

\subsection{A fundamental theorem of calculus}

We are ready for the main result.

\begin{theorem}
[Fundamental Theorem]\label{Main}Under the preceding Hypothesis,%
\[
\iint\nolimits_{f^{-1}[[a,b]]}D_{\mathbf{T}}^{2}f\,dA=\sigma\cdot
2\pi(b-a)\text{ .}
\]

\end{theorem}

\begin{proof}
Using (\ref{Eq Resolution}), Lemma \ref{Lemma Curvature of Level Curves}, and
(\ref{Eq Integral-of-Curvature}), we calculate as follows:%
\begin{equation}
\iint\nolimits_{f^{-1}[[a,b]]}D_{\mathbf{T}}^{2}f\,dA=\int_{a}^{b}\int%
_{f^{-1}[t]}\frac{D_{\mathbf{T}}^{2}f}{|\nabla f|}ds\,dt=\int_{a}^{b}%
\int_{f^{-1}[t]}\kappa ds\,dt=(2\pi\sigma)(b-a)\,\text{,}\nonumber
\end{equation}
establishing the claimed formula.
\end{proof}

We note an immediate consequence concerning the integral of $D_{\mathbf{N}%
}^{2}f$.

\begin{corollary}
Under the preceding assumptions,%
\[
\iint\nolimits_{f^{-1}[[a,b]]}D_{\mathbf{N}}^{2}f\,dA=\int_{f^{-1}%
[b]}\left\vert \nabla f\right\vert ds-\int_{f^{-1}[a]}\left\vert \nabla
f\right\vert ds-2\pi(b-a)\sigma\text{.}%
\]
If, in addition, $f$ is harmonic, then%
\[
\iint\nolimits_{f^{-1}[[a,b]]}D_{\mathbf{N}}^{2}f\,dA=-2\pi(b-a)\sigma\text{.}%
\]

\end{corollary}

We provide some hints for the proof and leave the details to the reader. Note
that the number $Q_{p}(\mathbf{i},\mathbf{i})+Q_{p}(\mathbf{j},\mathbf{j})$,
i.e., the trace of the Hessian form $Q_{p}$, equals the Laplacian $\Delta
f(p)$, and that the trace of a symmetric bilinear form is invariant under an
orthogonal change of coordinates. Hence,%
\[
D_{\mathbf{T}}^{2}f(p)+D_{\mathbf{N}}^{2}f(p)=Q_{p}(\mathbf{T},\mathbf{T}%
)+Q_{p}(\mathbf{N},\mathbf{N})=\operatorname*{Tr}Q_{p}=\Delta f(p)\text{.}%
\]
As $\Delta f=\operatorname{div}(\nabla f)$, we may apply Green's theorem to
the integral of $\Delta f$.

\section{Variations on the fundamental theorem}

We give three examples to illustrate some variations on the theme of Theorem
\ref{Main}.

\subsection{Examples}

In Examples 1 and 2, we will benefit from using complex-valued variables to
express real functions. In both examples, we adopt the following notations.

\smallskip

\noindent\textbf{Notation}. We name the Cartesian and polar coordinates of the
complex variables $z$ and $w$ as follows: $z=x+iy=\rho e^{i\varphi}$ and
$w=u+iv=re^{i\theta}$. We let $\overline{D}(0;t)$ denote the closed disc
$\{w:|w|\leq t\}$.

\smallskip

\noindent\textbf{Example 1}. Let $f(x,y)=|z+1|/|z-1|$, which is $C^{2}$ on $%
%TCIMACRO{\U{2102} }%
%BeginExpansion
\mathbb{C}
%EndExpansion
\setminus\{-1,1\}$. Except for $f^{-1}[0]$ (the singleton $\{-1\}$)\ and
$f^{-1}[1]$ (the $y$-axis), the level curves of $f$ are the so-called
\textit{Apollonius' circles} with foci $\pm1$. On the left\ half-plane, the
level circles are oriented counterclockwise according to our earlier
stipulation and $\sigma=1$ as a result. Exploiting the simple relation between
curvature of a circle and its radius, we find that $\kappa(z)=-2x/\left\vert
z^{2}-1\right\vert $.

We are to integrate $D_{\mathbf{T}}^{2}f$ over $f^{-1}[[0,1]]$, the entire
left half-plane. According to Theorem \ref{Main}, $\iint_{f^{-1}[[\epsilon
_{1},1-\epsilon_{2}]]}D_{\mathbf{T}}^{2}f\,dA=2\pi(1-\epsilon_{2}-\epsilon
_{1})$ for small positive $\epsilon_{1}$ and $\epsilon_{2}$. It follows that%
\[
\iint_{f^{-1}[[0,1]]}D_{\mathbf{T}}^{2}f\,dA=2\pi\,\text{.}
\]

We verify this claim by calculation. Introduce the auxiliary \textit{complex}
function $w(z)=(z+1)/(z-1)$, which maps $f^{-1}[[0,t]]$ \textit{one-to-one}
onto $\overline{D}(0;t)$. Taking advantage of Lemma
\ref{Lemma Curvature of Level Curves} and the fact that $|\nabla f|=\left\vert
w^{\prime}\right\vert $, we have%
\[
D_{\mathbf{T}}^{2}f\left(  z\right)  =\kappa(z)|\nabla f(z)|=\kappa
(z)\left\vert w^{\prime}(z)\right\vert \text{ .}%
\]
Change the variables of integration from $(x,y)$ to $(u,v)$ at the cost of
dividing the integrand by $\det\left[  \partial(u,v)/\partial(x,y)\right]
=|w^{\prime}|^{2}$, we obtain%
\[
\iint_{\operatorname{Re}z\leq0}D_{\mathbf{T}}^{2}f\,dA_{z}=\iint_{|w|\leq
1}\frac{D_{\mathbf{T}}^{2}f}{|w^{\prime}|^{2}}\,dA_{w}=\iint_{|w|\leq1}%
\frac{\kappa}{|w^{\prime}|}\,dA_{w}%
\]
(where the subscripts distinguish the area elements in the $z$-plane and
$w$-plane). Now,%
\[
\frac{\kappa}{|w^{\prime}|}=-\frac{2x}{|z^{2}-1|}\frac{|z-1|^{2}}%
{2}=-x\left\vert \frac{z-1}{z+1}\right\vert =-\frac{x}{|w|}=-\frac{x}{r}\text{
.}%
\]
Inverting the function $z\mapsto w(z)$ allows $x$ to be expressed in terms of
$w$:%
\[
x=\operatorname{Re}z=\operatorname{Re}\frac{w+1}{w-1}=\frac{|w|^{2}%
-1}{|w-1|^{2}}=\frac{r^{2}-1}{r^{2}-2r\cos\theta+1}\text{ .}%
\]
Finally, we have%
\[
\iint_{f^{-1}[[0,1]]}D_{\mathbf{T}}^{2}f\,dA=\int_{0}^{2\pi}\int_{0}^{1}%
\frac{1-r^{2}}{r^{2}-2r\cos\theta+1}drd\theta=2\pi\text{ ,}%
\]
as modern computing technology capable of symbolic integration can verify.

\smallskip

\noindent\textbf{Example 2}. Let $f(x,y)=|z^{2}-1|$, which is $C^{2}$ on $%
%TCIMACRO{\U{2102} }%
%BeginExpansion
\mathbb{C}
%EndExpansion
\setminus\{-1,1\}$ and has a saddle point at the origin. Let $C_{t}$ denote
$f^{-1}[t]$. For each $t>0$, the curve $C_{t}$ is a \textit{Cassini's
oval}\ with foci $\pm1$, which is the locus of points whose distances to
$\pm1$ have a fixed product $t$. For $t<1$, $C_{t}$ has two components, both
oriented counterclockwise by our stipulation; for $t=1$, $C_{t}$ has a single
self-intersection at the origin and is well known as a \textit{Bernoulli's
lemniscate}.

Let $\Omega=f^{-1}[(0,1)]$. In the spirit of\ Theorem \ref{Main}, we expect
that%
\[
\iint_{\overline{\Omega}}D_{\mathbf{T}}^{2}f\,dA=2\times2\pi\times
(1-0)=4\pi\text{ ,}
\]
because, for $t\in(0,1)$, $\int_{C_{t}}\kappa ds=2\times2\pi$, $C_{t}$ being
the union of \textit{two} positively-oriented simple closed curves. Let's
verify this conclusion by purely computational means independent of any
results established herein.

Define the auxiliary function $w(z)=z^{2}-1$, which is \textit{one-to-one} on
$R:=\left\{  \rho e^{i\varphi}:\rho\geq0;\,\varphi\in\left(  -\pi
/2,\pi/2\right]  \right\}  $ and maps $\overline{\Omega}\cap R$ \textit{onto}
$\overline{D}(0;1)$.

Using the Hessian form $Q$, we find%
\[
D_{\mathbf{T}}^{2}f\left(  z\right)  =2\frac{\left\vert z\right\vert
^{4}+(x^{2}-y^{2})}{\left\vert z\right\vert ^{2}|z^{2}-1|}=\frac{2\rho
^{4}+2\rho^{2}\cos2\alpha}{\rho^{2}r}\text{ .}
\]
Note that $r^{2}=|z^{2}-1|^{2}=\rho^{4}-2\rho^{2}\cos2\alpha+1$, enabling us
to rewrite $D_{\mathbf{T}}^{2}f\left(  z\right)  $:%
\[
D_{\mathbf{T}}^{2}f\left(  z\right)  =\frac{3\rho^{4}+1-r^{2}}{\rho^{2}%
r}\text{ .}
\]
Finally, we change variables from $(x,y)$ to $(u,v)$ and apply the relations
$\rho^{2}=|w+1|$ and $|w^{\prime}|=2\rho$ to obtain%
\begin{align*}
\iint\nolimits_{\overline{\Omega}\cap R}D_{\mathbf{T}}^{2}f\,dA_{z}  &
=\iint\nolimits_{\overline{D}(0;1)}\frac{D_{\mathbf{T}}^{2}f}{|w^{\prime}%
|^{2}}dA_{w}=\iint\nolimits_{\overline{D}(0;1)}\frac{3|w+1|^{2}+1-r^{2}%
}{4|w+1|^{2}r}dA_{w}\\
&  =\frac{1}{2}\int_{0}^{2\pi}\left(  \int_{0}^{1}1+\frac{r\cos\theta+1}%
{r^{2}+2r\cos\theta+1}dr\right)  d\theta\text{ .}%
\end{align*}
As, $\iint_{\overline{\Omega}}D_{\mathbf{T}}^{2}f\,dA=2\iint%
\nolimits_{\overline{\Omega}\cap R}D_{\mathbf{T}}^{2}f\,dA_{z}$, it suffices
to verify that%
\[
\int_{0}^{2\pi}\left(  \int_{0}^{1}\frac{r\cos\theta+1}{r^{2}+2r\cos\theta
+1}dr\right)  d\theta=2\pi\text{ ,}
\]
which our computing technology can attest to.

\smallskip

\noindent\textbf{Example 3}. Let $f(x,y)=[x^{2}+(y-1)^{2}-4][x^{2}%
+(y+1)^{2}-4]$. The exercise of deducing the following claims from Theorem
\ref{Main} is left to the reader.%
\[
\iint\nolimits_{f^{-1}[[-8,0]]}D_{\mathbf{T}}^{2}f\,dA=32\pi\,\text{; }%
\iint\nolimits_{f^{-1}[[0,8]]}D_{\mathbf{T}}^{2}f\,dA=0\,\text{; }%
\iint\nolimits_{f^{-1}[[1,20]]}D_{\mathbf{T}}^{2}f\,dA=22\pi\text{\thinspace.}%
\]

\subsection{Lessons from the examples}

With the preceding examples in mind, we discuss ways in which Theorem
\ref{Main} can be adapted. For convenience, we introduce two terms. If $p$ is
a critical point of $f$, then $f(p)$ is a \textit{critical value} of $f$ and
$f^{-1}[f(p)]$ is a \textit{critical level}.

Assume that $f$ has at most finitely many critical points in $\Omega$. We
further impose the condition (sufficient for our purpose) that any critical
level in $f^{-1}[[a,b]]$ has plane measure 0; a \textquotedblleft naturally
occurring\textquotedblright\ function easily meets this condition, as a
critical level is often a finite union of rectifiable arcs.

Under these assumptions, we outline three cases.

\noindent\textbf{Case 1. }$f^{-1}[[a,b]]$ is compact, free of critical points,
but disconnected. (Cf. Examples 2--3.) Theorem \ref{Main} can be applied to
each component $K_{j}$ of $f^{-1}[[a,b]]$, resulting in%
\begin{equation}
\iint\nolimits_{f^{-1}[[a,b]]}D_{\mathbf{T}}^{2}f\,dA=%
%TCIMACRO{\tsum \nolimits_{i}}%
%BeginExpansion
{\textstyle\sum\nolimits_{i}}
%EndExpansion
2\pi\sigma(K_{i})\cdot(b-a)\text{ .} \label{Eq Disconnected Domain}%
\end{equation}

\noindent\textbf{Case 2. }$f^{-1}[[a,b]]$ is compact and $[a,b]$ contains a
critical value. (Cf. Examples 2--3.) For simplicity but without loss of
generality, assume that $[a,b]$ contains exactly one critical value $c$.
Either $c$ is an endpoint or $c\in(a,b)$. The latter case can be reduced to
the former, as $[a,b]=[a,c]\cup\lbrack c,b]$. Suppose that $c=b$; the case
$c=a$ is similar. Then%
\[
\iint\nolimits_{f^{-1}[[a,b]]}D_{\mathbf{T}}^{2}f\,dA=\lim_{\epsilon
\rightarrow0}\iint\nolimits_{f^{-1}[[a,b-\epsilon]]}D_{\mathbf{T}}%
^{2}f\,dA\text{ .}
\]
As the integral on the right is\ governed by (\ref{Eq Disconnected Domain}),
the limit exists. (Our assumption that a critical level has plane measure $0$
is needed for the above equality. Also note that, although $D_{\mathbf{T}}%
^{2}f$ is undefined at the critical points on the critical level $f^{-1}[b]$,
its integral over $f^{-1}[[a,b]]$ after all exists.)

\noindent\textbf{Case 3. }$[a,b]$ contains no critical values but, for
finitely many $c\in\lbrack a,b]$, $f^{-1}[c]$ is not compact. (Cf. Example 1.)
As in Case 2, assume that there is exactly one such $c$ and $c=b$. The same
limit argument in Case 2 shows that $D_{\mathbf{T}}^{2}f$ (defined everywhere
on $f^{-1}[[a,b]]$) is integrable on $f^{-1}[[a,b]]$ and given by
(\ref{Eq Disconnected Domain}).

\end{document}